\documentclass{article}
\usepackage[centertags]{amsmath}

\begin{document}

\title{The geometrical origins of some distributions and the complete concentration of measure phenomenon for mean-values of functionals}

\author{ Cheng-shi Liu\\Department of Mathematics\\Northeast Petroleum University\\Daqing 163318, China
\\Email:chengshiliu-68@126.com}

\maketitle

 We derive out naturally some important distributions such as high order normal distributions and high order exponent distributions and the Gamma distribution from a geometrical way. Further, we obtain the exact mean-values of integral form functionals in the balls of continuous functions space with $p-$norm, and show the complete concentration of measure phenomenon which means that a functional takes its average on a ball with probability 1, from which we have nonlinear exchange formula of expectation.



\section{Introduction}

Following the works of  G\^{a}teaux,  L\'{e}vy and Wiener et al[1-9], we study the mean-values of some functionals in the ball of continuous functions space under $p-$norm. The computation of the mean-values of functionals is still an important problem for mathematicians and physicians. Some infinite-dimensional probability problems lead to naturally the averages of functionals. For example, if we choose randomly a curve $x(t)$ in $M=\{x|0\leq x'(t)\leq1\}$, what is the mean-value of its arc length $l(x)=\int_0^1 \sqrt{1+(x'(t))^2}\mathrm{d}t$ on $M$? what is the mean-value of its area $A(x)=\int_0^1 x(t)\mathrm{d}t$ on $M$? We can get the mean-value $E(l(x))=\frac{\sqrt2}{2}+\frac{1}{2}\ln(1+\sqrt2)$, and $E(A(x))=\frac{1}{4}$.  Although these problems belong to essentially the infinite-dimensional integrations, they are different with the functional integrations such as Feynman integral and Wiener integral in quantum physics and stochastic process[10-13]. In general, the mean-value of functional $f$
on $M$ can be defined as follows,
\begin{equation}
Ef=\frac{\int_{M} f(x)D(x)}{\int_{M}D(x)},
\end{equation}
where $D(x)$ represents formally  the volume
element of $M$. However, there exists a huge gap between the definitions and computations of functional integrations. In particular, for the normed (or more general, topology) linear space $V$, the integration on it is defined by cylinder measure depending on the dual space $V^*$ ([11,13]). For example, for the space $V=L^p[a,b]$, the construction of measure on it is so complicated that one almost cannot use it to compute the concrete integrations of some simple functionals  such as
$f[x]=\int_a^bx(t)\mathrm{d}t$, where we require $x(t)\in M=\{x(t)|0\leq x(t)\leq 1, a\leq t\leq b, x(t)\in L[a,b]\}$.

The computation of mean-value of a functional on a set $M$ depends on not only the form of functional, but also the structure of the set $M$. Abstractly, a functional
$f(x)$ is a function of $x$ where $x$ is an element in an
infinite-dimensional space such as $C[0,1]$.  In
general, there are two basic ways to construct functionals. One
way is to use the values of $x(t)$ on some points
$t_1,\cdots,t_m$ such that $f(x)=g(x(t_1),\cdots,x(t_m))$, where $g$ is a usual function in $R^n$. Essentially, all the first kind of functionals are finite dimensional functions. Another way is to use
integral
\begin{equation}
f(x)=\int_{I_1}\cdots\int_{I_m}g(x(t_1),\cdots,x(t_m))\mathrm{d}t_1\cdots\mathrm{d}t_m,
\end{equation}
where $I_1,\cdots,I_m$ are subsets of the interval $[0,1]$. In
general, we take every $I_k$ to be a subinterval.  For example,
\begin{equation*}
f(x)=\int_0^{0.3}x^2(t)\mathrm{d}t+\sin(\int_{0.1}^{0.2}x(t)\mathrm{d}t+\int_0^{0.8}\cos x(t)\mathrm{d}t).
\end{equation*}
Therefore, there are two kinds of basic elements $x(t_i)(i=1,\cdots,n)$ and
$\int_{I_j}g(x(t))\mathrm{d}t$ $(j=1,\cdots,m)$ such that those useful functionals
can be constructed in terms of them by addition, subtraction, multiplication, division and composition. We also call these functionals
the elementary functionals.

 For the first kind of functionals, their averages are just reduced to the usual finite-dimensional
integrals. Thus we only consider the integrals of the second
kind of functionals $f(x)$. If the domain of the functional $f$ is $M$,  in general, the volume $\int_{M}D(x)$ of $M$ is zero or
infinity, and the infinite-dimensional integral $\int_{M}
f(x)D(x)$ is always respectively zero or infinity (or finite sometimes). However, the mean-value $
Ef=\frac{\int_{M} f(x)D(x)}{\int_{M}D(x)}$
perhaps is finite. Therefore we  need a reasonable
 definition of the mean-value of the functional. The standard approach is to use a limit procedure.
 For example, if we take the set $M=\{x|a\leq x(t)\leq b, x(t)\in C[0,1]\}$ and the functional $f(x)=\int_0^1g(x(t))\mathrm{d}t$,
 then we define the mean-value of $f$ as
\begin{equation}
Ef=\lim_{n\rightarrow \infty}\frac{\int_a^b\cdots\int_a^b \frac{1}{n}\sum_{k=1}^ng(x_k)\mathrm{d}x_1\cdots\mathrm{d}x_n}{\int_a^b\cdots\int_a^b \mathrm{d}x_1\cdots\mathrm{d}x_n}
\end{equation}
where $x_k=x(\frac{k}{n})$. If the limitation exists, we call it the mean-value of the functional $f$.  In addition, we must emphasize that for different function spaces, we need different limitation procedures.

In the paper,  we mainly study the mean-values of the most important functionals $f$ of form (2) in the ball $M$ of the continuous functions space under $p$-norm.  Our method is to consider the functional as an infinite-dimensional
random variable and then to compute the expectation  and the variances on $M$. For the purpose, we first obtain the probability distribution of the coordinate $x(t)$  in $M$ for fixed $t$. These distributions have the forms of high order normal and exponent distributions by which we can also derive naturally the important Gamma distribution in statistics. Then we show that the variances are zeros which indicate that these functionals satisfy the property of the complete concentration of measure which means that a functional takes its mean-value on an infinite-dimensional ball with probability 1. Corresponding, we give the nonlinear exchange formulas for averages of functionals. The routine concentration of measure is described by some inequalities such as Levy lemma[14-16], which is different with the complete concentration of measure which is shown by some equalities.

\section{The geometrical origins of some distributions}

Consider the continuous functions space $C[0,1]$ and define some norms such as $||x||_0=\max_{t\in[0,1]}|x(t)|$, and $||x||_p=(\int_0^1|x(t)|^p\mathrm{d}t)^{\frac{1}{p}}$ for $p\geq1$. For $p=\frac{p_0}{q_0}$ where $p_0$ is even and $(p_0,q_0)=1$, we consider the whole ball
$M=\{x|||x||_p\leq R^p\}$, while for $p\geq1$ is a general real number or  specially $p=\frac{p_0}{q_0}$ where  $p_0$ is odd and $(p_0,q_0)=1$, we only consider the "first quadrant" of $M$, that is
$M^+=\{x|x\geq0,||x||_p\leq R^p\}$.

The following lemma is important.

\textbf{Lemma}([17]). The following generalized Dirichlet formula holds
\begin{equation}
\int\cdots\int_{B^+}x_1^{p_1-1}x_2^{p_2-1}\cdots x_n^{p_n-1}\mathrm{d}x_1\cdots\mathrm{d}x_n=
\frac{1}{2^n}\frac{\Gamma(\frac{p_1}{2})\cdots \Gamma(\frac{p_n}{2})}{\Gamma(1+\frac{p_1+\cdots+p_n}{2})},
\end{equation}
where $p_i>0$ for $i=1,\cdots,m$ and $B^+=\{(x_1,\cdots,x_n)|x_1^2+\cdots+x_n^2\leq1, x_k>0, k=1,\cdots, n\}$.

Next we give the following  results.

\textbf{Theorem 1}. For the ball $M=\{x|\int_0^1x^p(t)\mathrm{d}t\leq R^p\}$ in $C[0,1]$ when $P=\frac{p_0}{q_0}$ where $p_0$ is even and $(p_0,q_0)=1$, the density of $x(t)$ as a random variable for fixed $t$ is given  by
\begin{equation}
\rho(x)=\frac{1}{2R\Gamma(\frac{1}{p})p^{\frac{1}{p}-1}}\mathrm{e}^{-\frac{x^p}{pR^p}}, x\in(-\infty, +\infty).
\end{equation}

\textbf{Proof}. Firstly, by discretization, we have $M_n=\{(x_1,\cdots,x_n)|x_1^p+\cdots+x_n^p\leq nR^p\}$ where
$x_k=x(\frac{k}{n})$. We direct compute the density $\rho_n(x)$ of $x_1$ as a random variable in the ball $M_n$. We have from the above lemma
\begin{equation*}
\rho_n(x_1)=\frac{\int\cdots\int_{M_n}\mathrm{d}x_2\cdots\mathrm{d}x_n}{\int_{M_n}\mathrm{d}V}
\end{equation*}
\begin{equation}
=\frac{p\Gamma(1+\frac{n}{p})}{2\Gamma(\frac{1}{p})\Gamma(1+\frac{n-1}{p})(nR^p-x^2_1)^{\frac{1}{p}}}
(1-\frac{x_1^2}{nR^p})^{\frac{n}{p}}.
\end{equation}
Taking the limitation of $n$ approaching to $+\infty$, and using the Stirling's asymptotic formula of Gamma function,  we get for $t$ fixed the density of $x(t)$ in $M$,
\begin{equation*}
\rho(x)=\lim_{n\rightarrow+\infty}\frac{p\sqrt{2\pi}\sqrt{\frac{n}{p}}
(\frac{n}{p})^{\frac{n}{p}}\mathrm{e}^{-\frac{n}{p}}}
{2\Gamma(\frac{1}{p})\sqrt{2\pi}\sqrt{\frac{n-1}{p}}
(\frac{n-1}{p})^{\frac{n-1}{p}}\mathrm{e}^{-\frac{n-1}{p}}(nR^p-x^2)^{\frac{1}{p}}}
(1-\frac{x^2}{nR^p})^{\frac{n}{p}}
\end{equation*}
\begin{equation}
=\lim_{n\rightarrow+\infty}\frac{1}{2Rp^{\frac{1}{p}-1}\Gamma(\frac{1}{p})\mathrm{e}^{\frac{1}{p}}}(\frac{n-1}{nR^p-x^p})^{\frac{1}{p}}
\sqrt{\frac{n}{n-1}}(\frac{n}{n-1})^{\frac{n}{p}}
(1-\frac{x^2}{nR^2})^{\frac{n}{p}}
\end{equation}
\begin{equation}
=\frac{1}{2Rp^{\frac{1}{p}-1}\Gamma(\frac{1}{p})}\mathrm{e}^{-\frac{x^p}{pR^p}}.
\end{equation}
  The proof is completed.

\textbf{Theorem 2}. For  the "first quadrant" $M^+=\{x|x(t)\geq0, \int_0^1x^p(t)\mathrm{d}t\leq R^p\}$ of $M$ when $p$ is a general real number and $p\geq 1$ or specially $p=\frac{p_0}{q_0}$ where $p_0$ is odd and $(p_0,q_0)=1$, the density of $x(t)$ as a random variable on $M^+$ for fixed $t$ is given by
\begin{equation}
\rho(x)=\frac{1}{Rp^{\frac{1}{p}-1}\Gamma(\frac{1}{p})}\mathrm{e}^{-\frac{x^p}{pR^p}}, x\in[0, +\infty).
\end{equation}

\textbf{Proof}.  Firstly, by discretization, we have $M^+_n=\{(x_1,\cdots,x_n)|x_1^p+\cdots+x_n^p\leq nR^p, x_k\geq 0, k=1,\cdots,n\}$ where
$x_k=x(\frac{k}{n})$. We direct compute the density $\rho_n(x)$ of $x_1$ as a random variable in the "first quadrant" $M^+_n$. We have from the lemma
\begin{equation*}
\rho_n(x_1)=\frac{\int\cdots\int_{M^+_n}\mathrm{d}x_2\cdots\mathrm{d}x_n}{\int_{M^+_n}\mathrm{d}V}
\end{equation*}
\begin{equation}
=\frac{p\Gamma(1+\frac{n}{p})}{\Gamma(\frac{1}{p})\Gamma(1+\frac{n-1}{p})(nR^p-x^2_1)^{\frac{1}{p}}}
(1-\frac{x_1^2}{nR^p})^{\frac{n}{p}}.
\end{equation}
Taking the limitation of $n$ approaching to $+\infty$, and using the Stirling's asymptotic formula of Gamma function,  we get for $t$ fixed the density of $x(t)$ in $M$,
\begin{equation*}
\rho(x)=\lim_{n\rightarrow+\infty}\frac{p\sqrt{2\pi}\sqrt{\frac{n}{p}}
(\frac{n}{p})^{\frac{n}{p}}\mathrm{e}^{-\frac{n}{p}}}
{\Gamma(\frac{1}{p})\sqrt{2\pi}\sqrt{\frac{n-1}{p}}
(\frac{n-1}{p})^{\frac{n-1}{p}}\mathrm{e}^{-\frac{n-1}{p}}(nR^p-x^2)^{\frac{1}{p}}}
(1-\frac{x^2}{nR^p})^{\frac{n}{p}}
\end{equation*}
\begin{equation}
=\lim_{n\rightarrow+\infty}\frac{1}{Rp^{\frac{1}{p}-1}\Gamma(\frac{1}{p})\mathrm{e}^{\frac{1}{p}}}(\frac{n-1}{nR^p-x^p})^{\frac{1}{p}}
\sqrt{\frac{n}{n-1}}(\frac{n}{n-1})^{\frac{n}{p}}
(1-\frac{x^2}{nR^2})^{\frac{n}{p}}
\end{equation}
\begin{equation}
=\frac{1}{Rp^{\frac{1}{p}-1}\Gamma(\frac{1}{p})}\mathrm{e}^{-\frac{x^p}{pR^p}}.
\end{equation}
  The proof is completed.

\textbf{Remark 1}. When $p$ is even, the density looks like a normal distribution, and thus we call it high order normal distribution or normal-like distribution.  If $R=1$, a simple form is
\begin{equation}
\rho(x)=\frac{1}{2\Gamma(\frac{1}{p})p^{\frac{1}{p}-1}}\mathrm{e}^{-\frac{x^p}{p}}, x\in(-\infty, +\infty).
\end{equation}
Further, for example, if $p=2$, we get the standard normal distribution
\begin{equation}
\rho(x)=\frac{1}{\sqrt{2\pi}}\mathrm{e}^{-\frac{x^2}{2}}, x\in(-\infty, +\infty),
\end{equation}
which gives the G\^{a}teaux and L\'{e}vy's result [1].
If $p=4$, we get a  4-order normal distribution
\begin{equation}
\rho(x)=\frac{\sqrt2}{\Gamma(\frac{1}{4})}\mathrm{e}^{-\frac{x^4}{4}}, x\in(-\infty, +\infty).
\end{equation}

\textbf{Remark 2}. When $p$ is odd, the density looks like an exponent distribution, and thus we call it high order exponent distribution or exponent-like distribution.  For example, if $p=1$ and $R=\lambda^{-\frac{1}{p}} $, we get the usual exponent distribution
\begin{equation}
\rho(x)=\frac{1}{\lambda}\mathrm{e}^{-\lambda x}, x\in[0, +\infty).
\end{equation}
If $p=3$ and $R=(3\lambda)^{-\frac{1}{p}} $, we get the 3-order exponent distribution
\begin{equation}
\rho(x)=\frac{3\lambda^{\frac{1}{3}}}{\Gamma(\frac{1}{3})}\mathrm{e}^{-\lambda x^3}, x\in[0, +\infty).
\end{equation}

\textbf{Remark 3}. By a simple transformation, we can obtain the famous Gamma distribution in statistics. Indeed, we take a transformation
\begin{equation}
Z=\frac{x^p}{p\beta},
\end{equation}
then the density of $Z$ is just
\begin{equation}
\rho(z)=\frac{\beta^{\frac{1}{p}}}{\Gamma(\frac{1}{p})}z^{\frac{1}{p}-1}\mathrm{e}^{-\beta y}.
\end{equation}
Further, taking $\alpha=\frac{1}{p}$ gives the Gamma distribution
\begin{equation}
\rho(z,\alpha,\beta)=\frac{\beta^{\alpha}}{\Gamma(\alpha)}z^{\alpha-1}\mathrm{e}^{-\beta y}.
\end{equation}
 This is a geometrical origin of the Gamma distribution. We can see that this is a rather natural way to derive the Gamma distribution.

\section{The exact mean-values of functionals and the complete concentration of measure phenomenon}

\textbf{Theorem 3}. Take the ball $M=\{x|\int_0^1x^p(t)\mathrm{d}t\leq R^p\}$ in $C[0,1]$ with norm $||x||_p$ when  $p=\frac{p_0}{q_0}$ where $p_0$ is even and $(p_0,q_0)=1$. Then the mean-value and invariance on $M$ of the functional
\begin{equation}
Y=f(x)=\int_0^1\cdots \int_0^1g(x(t_1),\cdots, x(t_m))\mathrm{d}t_1\cdots\mathrm{d}t_m,
\end{equation}
satisfy
\begin{equation}
EY=\frac{1}{2^m\Gamma^m(\frac{1}{p})p^{\frac{m}{p}-m}}\int_{-\infty}^{+\infty}\cdots\int_{-\infty}^{+\infty}
g(Rx_1,\cdots,Rx_m)\mathrm{e}^{-\frac{x_1^p}{p}-\cdots-\frac{x_m^p}{p}}\mathrm{d}x_1\cdots\mathrm{d}x_m,
\end{equation}
\begin{equation}
DY=0.
\end{equation}
Further, for the functionals $Y_1,\cdots, Y_m$ with the form (2) and a general function $h$ of $m$ variables,  we have the nonlinear exchange formula for the mean-value of $h(Y_1,\cdots,Y_m)$ on $M$,
\begin{equation}
Eh(Y_1,\cdots,Y_m)=h(EY_1,\cdots,EY_m).
\end{equation}

\textbf{Proof}. We only need to prove it in the case of $m=1$. By the linearity of the expectation $E$, we have
\begin{equation}
EY=\int_0^1E(g(x(t))\mathrm{d}t=E(g(x(t))
=\frac{1}{2R\Gamma(\frac{1}{p})p^{\frac{1}{p}-1}}\int_{-\infty}^{+\infty}g(Rx)\mathrm{e}^{-\frac{x^p}{pR^p}}\mathrm{d}x.
\end{equation}
Then replacing $x$ by $px$, we get the result. Similarly, by the two-dimensional measure of the set $\{(t,s)|t=s\}$ as a subset of $[0,1]^2$ being zero, we have,
\begin{equation}
E(Y^2)=\int_0^1\int_0^1E(g(x(t)g(x(s))\mathrm{d}t\mathrm{d}s=E(g(x(t))E(g(x(s))=E^2(Y),
\end{equation}
and then $DY=0$ by which the nonlinear exchange formula holds naturally. The proof is completed.

Similarly, we have the following theorem. 

\textbf{Theorem 4}. Take the "first quadrant" $M^+=\{x|x(t)\geq0, \int_0^1x^p(t)\mathrm{d}t\leq R^p\}$ of $M$ when $p=\frac{p_0}{q_0}$ where $p_0$ is odd and $(p_0,q_0)=1$. Then the mean value and invariance on $M$ of the functional
\begin{equation}
Y=f(x)=\int_0^1\cdots \int_0^1g(x(t_1),\cdots, x(t_m))\mathrm{d}t_1\cdots\mathrm{d}t_m,
\end{equation}
satisfy
\begin{equation}
EY=\frac{1}{\Gamma^m(\frac{1}{p})p^{\frac{m}{p}-m}}\int_{0}^{+\infty}\cdots\int_{0}^{+\infty}
g(Rx_1,\cdots,Rx_m)\mathrm{e}^{-\frac{x_1^p}{p}-\cdots-\frac{x_m^p}{p}}\mathrm{d}x_1\cdots\mathrm{d}x_m,
\end{equation}
\begin{equation}
DY=0.
\end{equation}
Further, for the functionals $Y_1,\cdots, Y_m$ with the form (2) and a general function $h$ of $m$ variables,  we have the nonlinear exchange formula for the average of $h(Y_1,\cdots,Y_m)$ on $M$,
\begin{equation}
Eh(Y_1,\cdots,Y_m)=h(EY_1,\cdots,EY_m).
\end{equation}

\textbf{Remark 4}. The theorems 3 and 4 show the complete concentration of measure phenomenon which means that a functional takes its mean-value on an infinite-dimensional ball with probability 1.

\textbf{Remark 5}. In infinite-dimensional sphere, we have the routine concentration of measure. Indeed, for the sphere $M(R)=\{x|\int_0^1x^p(t)\mathrm{d}t\leq R^p\}$ ($p=\frac{p_0}{q_0}$ where $p_0$ is even and $(p_0,q_0)=1$) or $M(R)=\{x|\int_0^1x^p(t)\mathrm{d}t\leq R^p, x(t)\geq0\}$ ($p=\frac{p_0}{q_0}$ where $p_0$ is odd and $(p_0,q_0)=1$) in $C[0,1]$ with $p-$norm, we can easily see that the measure concentrates on the surface of sphere $M$ since $\frac{V(M(R))}{V(M(r))}=\lim_{n\rightarrow+\infty}(\frac{r}{R})^n=0$ for $r<R$. Therefore, the mean-values of the functional in the set $M$ (or $M^+$) can be reduced to and equal to the mean-value on the surface of $M$(or $M^+$).

\textbf{Remark 6}. For the functional
\begin{equation}
Y=f(x)=\int_{I_1}\cdots\int_{I_m}g(x(t_1),\cdots,x(t_m),t_1,\cdots,t_m)\mathrm{d}t_1\cdots\mathrm{d}t_m,
\end{equation}
in general, the mean-value on $M$ or $M^+$ will be respectively
\begin{equation*}
E(Y)=\frac{1}{2^m\Gamma^m(\frac{1}{p})p^{\frac{m}{p}-m}}\int_{I_1}\cdots\int_{I_m}\int_{-\infty}^{+\infty}\cdots\int_{-\infty}^{+\infty}
g(Rx_1,\cdots,Rx_m,t_1,\cdots,t_m)
\end{equation*}
\begin{equation}
\times\mathrm{e}^{-\frac{x_1^p}{p}-\cdots-\frac{x_m^p}{p}}
\mathrm{d}x_1\cdots\mathrm{d}x_m\mathrm{d}t_1\cdots\mathrm{d}t_m,
\end{equation}
or
\begin{equation*}
E(Y)=\frac{1}{\Gamma^m(\frac{1}{p})p^{\frac{m}{p}-m}}\int_{I_1}\cdots\int_{I_m}\int_{0}^{+\infty}\cdots\int_{0}^{+\infty}
g(Rx_1,\cdots,Rx_m,t_1,\cdots,t_m)
\end{equation*}
\begin{equation}
\mathrm{e}^{-\frac{x_1^p}{p}-\cdots-\frac{x_m^p}{p}}
\mathrm{d}x_1\cdots\mathrm{d}x_m\mathrm{d}t_1\cdots\mathrm{d}t_m.
\end{equation}
But in general $DY\neq0$ and then the complete concentration of measure don't hold.

\end{document}